\documentclass[11pt]{amsart}
\usepackage{latexsym}
\usepackage{amsfonts}

\usepackage{amsmath,amsthm,amssymb}

\title[Random quotients  of the modular group]{Random quotients  of the modular group are rigid and essentially
incompressible}

\author[I.~Kapovich]{Ilya Kapovich}

\address{\tt Department of Mathematics, University of Illinois at
  Urbana-Champaign, 1409 West Green Street, Urbana, IL 61801, USA
\newline http://www.math.uiuc.edu/\~{}kapovich/}
\email{\tt kapovich@math.uiuc.edu}

\curraddr{\tt Department of Mathematics, Johann Wolfgang Goethe
University, Robert Mayer Strasse 6-8, 60325 Frankfurt, Germany}

\author[P.~Schupp]{Paul Schupp}

\address{\tt Department of Mathematics, University of Illinois at
  Urbana-Champaign, 1409 West Green Street, Urbana, IL 61801, USA
\newline http://www.math.uiuc.edu/People/schupp.html}
\email{schupp@math.uiuc.edu}

\thanks{The authors were supported by the NSF
  grant DMS\#0404991. The first author also acknowledges the support of the Humboldt Foundation Research Fellowship.}

\newtheorem{theor}{Theorem}

\newtheorem{thm}{Theorem}[section]
\newtheorem{lem}[thm]{Lemma}
\newtheorem{cor}[thm]{Corollary}
\newtheorem{conj}[thm]{Conjecture}
\newtheorem{prop}[thm]{Proposition}
\theoremstyle{definition}
\newtheorem{defn}[thm]{Definition}
\newtheorem{notation}[thm]{Notation}
\newtheorem{conv}[thm]{Convention}
\newtheorem{rem}[thm]{Remark}

\newtheorem{prop-defn}[thm]{Proposition-Definition}

\begin{document}

\begin{abstract}
     We show that for any positive integer  $m\ge 1$,
  $m$-relator quotients of the modular group $M = PSL(2,\mathbb{Z})$
  generically satisfy a very strong Mostow-type \emph{isomorphism rigidity}.
  We also prove that such quotients
   are generically ``essentially  incompressible''. By this we mean that
  their ``absolute $T$-invariant'', measuring the smallest size of any
  possible finite presentation of the group,
  is  bounded below by a function which is
  almost linear  in terms of the length of the  given
  presentation.  We  compute the precise asymptotics of  the
  number $I_m(n)$ of \emph{isomorphism types} of $m$-relator quotients
  of $M$ where all the defining relators are cyclically reduced words
  of length $n$ in $M$.   We   obtain  other algebraic
  results and show that such quotients are complete, Hopfian, co-Hopfian,
  one-ended,  word-hyperbolic groups.
\end{abstract}


\maketitle


\section{Introduction}\label{intro}

The idea of genericity in Geometric Group Theory, understood as
the study of algebraic properties of random group-theoretic
objects, was introduced by Gromov~\cite{Grom,Grom1} when he
indicated that finitely presented groups are ``generically''
word-hyperbolic. This approach was made precise by
Ol'shanskii~\cite{Ol92}, Arzhantseva-Ol'shanskii~\cite{AO} and
Champetier~\cite{Ch94,Ch95}. Investigations centered around
genericity are now  an active and important research area (see,
for
example~\cite{AO,A1,A2,A3,Che98,Olliv,Olliv1,Olliv2,Olliv3,Z}).
One of the main reasons for studying genericity is that one can
use  the probabilistic method to discover  the existence of
objects with  new and interesting algebraic, algorithmic and
geometric properties. A  major example is Gromov's recent
construction~\cite{Grom2} of  a finitely presented group that is
not uniformly embeddable into a Hilbert space, which is  related
to  possible counter-examples   to the Novikov conjecture.

  There are many additional aspects of randomness
and genericity.  Work of the authors with   other colleagues,
~\cite{KMSS,KMSS1,KS,KSS,KSn},   introduced the notion of
\emph{generic-case complexity} for  decision problems.
It turns out that most classic group-theoretic decision problems,
such as the word, conjugacy, membership and the isomorphism problems,
 have \emph{provably}  low complexity on ``random'' inputs
even if  their  worst-case complexity is very high or even
unsolvable. This work subsequently led us to the discovery of
``isomorphism rigidity'' for generic  groups~\cite{KS,KSS}. The
famous Mostow Rigidity Theorem~\cite{Mos73} states  that if $M_1$
and  $M_2$ are complete connected hyperbolic manifolds of finite
volume and  dimension $n\ge 3$ then their fundamental groups are
isomorphic if and only if the manifolds themselves are isometric.
For a  group $G$ with a given finite generating set $A$ the
naturally associated geometric structure  is the Cayley graph
$\Gamma(G,A)$.  We thus say that a class of groups  equipped with
a specified finite generating set $A$  is \emph{rigid} if whenever
two groups from this class are isomorphic then their Cayley graphs
on  $A$  with the word metric are isometric.   Phenomena of this
type were known for various classes of Coxeter and Artin groups
(e.g \cite{Ros,PrSp,Bahls,BMMN,MW}).  In \cite{KSS} Kapovich,
Schupp and Shpilrain proved  the first  such theorem for a
``general'' class of groups by establishing  ``isomorphism
rigidity''  for generic one-relator groups.

   In most  previously explored contexts rigidity comes from a
careful analysis of some kind of particular structure. Thus  proofs of
quasi-isometric rigidity for non-uniform lattices in semi-simple Lie
groups hinge on the study of the structure of flats (see, for
example~\cite{Schw}).  Generic groups provide a conceptually new
source of group-theoretic rigidity where rigidity comes from
the properties of randomness itself.   Rigidity then opens
the way to proving results about  essential incompressibility
and the exact asymptotics of the number of isomorphism classes.

Before stating the main results we  introduce some definitions and
notation.

\begin{conv}[The Modular group]
 It is well-known that  the modular group $PSL(2,\mathbb Z)$ is isomorphic
to the free product
\[
M=\langle a,b| a^2=b^3=1\rangle.
\]
For the remainder of the paper we identify the modular group with $M$ and
use the free product structure.
We use  $\eta:M\to M$ to denote the automorphism of $M$ defined by
$\eta(a)=a$, $\eta(b)=b^{-1}$.

We  take  the group alphabet to be $A=\{a,b, b^{-1}\}\subseteq M$.
If  $w$ is a word in the alphabet $A$ then  $|w|$ denotes the length of $w$.
A  word $w$ in the alphabet $A$ is \emph{reduced} if it
does not contain any subwords of the form $aa,bb^{-1}, b^{-1}b,bb,
b^{-1}b^{-1}$. A  word $w$ is \emph{cyclically reduced} if
all cyclic permutations of $w$ are reduced. Note that if $|w|>1$
and $w$ is cyclically reduced then $|w|$ is necessarily even.
As for free groups, every   element $g\in M$ is represented by a unique reduced word $w$ in
the alphabet $A$  and we define  $|g|=|w|$.

If $G$ is a group and $R\subseteq G$, we denote by $\langle
\langle R\rangle\rangle$ the normal closure of $R$ in $G$, that
is, the smallest normal subgroup of $G$ containing $R$.

\end{conv}

\begin{notation}
We denote the set of all cyclically reduced words in the alphabet
$A$ by $\mathcal C$. If $\tau=(r_1,\dots, r_m)\in \mathcal C^m$ is
an $m$-tuple of cyclically reduced words, the  \emph{symmetrized
closure}, $R(\tau)$,  of  $\tau$ in $M$  is the set consisting of
all cyclic permutations of the elements of  $\tau$ and their
inverses. A set is \emph{symmetrized} if it is already equal to
its symmetrized closure.

If  $\tau=(r_1,\dots, r_m)\in M^m$  an $m$-tuple define
\[
|\tau|:=\max_{1 \le i \le m}  |r_i|
\]

 We set
\[
G_\tau:=M/\langle \langle \tau\rangle\rangle=\langle a,b\ |\
a^2=b^3=r_1=\dots=r_m=1\rangle.
\]

 If $m\ge 1$ is an integer then  $\mathcal T_m$
denotes the set of all $m$-tuples $(r_1,\dots, r_m)\in \mathcal
C^m$ such that $|r_1|=\dots=|r_m|$.
\end{notation}

  We next state the definition of ``genericity'' which we are using.

\begin{conv}
Let $m\ge 1$ and let $S\subseteq (A^*)^m$ be a set of $m$-tuples
of words in the alphabet $A$ and let $n\ge 0$. Then
\[
\gamma(n,S)=\#\{\tau\in S: |\tau|=n \}
\]
and
\[
\rho(n,S)=\#\{\tau\in S: |\tau|\le n\}.
\]
\end{conv}

\begin{defn}[Genericity]
Let $m\ge 1$ be an integer and let $X\subseteq M^m$ be a nonempty
subset. Let $S\subseteq X$.

We say that $S$ is \emph{exponentially generic} in $X$ if
\[
\lim_{n\to\infty} \frac{\rho(n,S)}{\rho(n,X)}=1,
\]
and the convergence is exponentially fast.
Similarly, we say that
$S$ is \emph{exponentially negligible} in $X$ if
\[
\lim_{n\to\infty} \frac{\rho(n,S)}{\rho(n,X)}=0,
\]
and the convergence is exponentially fast.
\end{defn}

Clearly, a subset of $X$ is exponentially generic if and only if its
complement in $X$ is exponentially negligible.
We will only be interested in the cases $X=\mathcal C^m$ and
$X=\mathcal T_m$ in this paper.

Our first main result about random quotients of the modular group
is:

\begin{theor}\label{A}[Isomorphism Rigidity]
  Let $m\ge 1$.  There exist an exponentially generic subset
  $Q_m$ of $\mathcal C^m$ and an exponentially generic subset
  $U_m=Q_m\cap \mathcal T_m$ of $\mathcal T_m$ such that the
  following hold:

\begin{enumerate}
\item There exists an algorithm that, given $\tau\in \mathcal
C^m$, decides, in time quartic  in  $|\tau|$, whether or not
$\tau\in Q_m$.

\item For every $\tau\in Q_m$ the group $G_\tau$ is
word-hyperbolic
  and one-ended. Moreover, the elements $a$ and $b$ respectively have orders $2$
  and $3$ in $G_\tau$.

\item For every $\tau\in Q_m$ the group $G_\tau$ is complete,
that is, the center of $G_\tau$ is trivial and $Out(G_\tau)=1$.

\item For any $\tau\in Q_m$ and \emph{any} finite symmetrized
  $S\subseteq \mathcal C$ satisfying $C'(1/8)$ we have $G_\tau\simeq
  M/\langle\langle S\rangle\rangle$ if and only if $R(\tau)=S$ or $R(\tau)=\eta(S)$ in $M$.
\item For $\tau=(r_1,\dots, r_m), \sigma=(s_1,\dots,s_m)\in Q_m$ we
  have $G_\tau\simeq G_\sigma$ if and only if there exist a reordering
  $\tau'=(r_1',\dots, r_m')$ of $\tau$ and $\epsilon\in\{0,1\}$ such
  that each $r_i'$ is a cyclic permutation of $\eta^\epsilon(s_i)$ or
  $\eta^\epsilon(s_i^{-1})$ for $i=1,\dots, m$.
\item If $\tau\in U_m, \sigma\in U_p$ with $|\tau|=|\sigma|$ but  $p>m$
  then $G_\tau\not\simeq G_\sigma$.
\end{enumerate}

\end{theor}

Part (5) of the above theorem says that for $\sigma, \tau\in Q_m$
the groups $G_\sigma$ and $G_\tau$ are isomorphic if and only if
their Cayley graphs with respect to the \emph{given} generating
set $a,b$ are isomorphic as labelled graphs by a graph
isomorphism preserves labels $a$ of edges and either preserves all
labels $b$ or inverts all of them.

\begin{theor}\label{B}[Homomorphism Rigidity]
  Let $m\ge 1$ be an integer.  Then the following hold:
\begin{enumerate}
\item Let $\tau\in Q_m$. Then for any homomorphism
  $\psi: M\to G_\tau$ exactly one of the following occurs:

(a) The image $\psi(M)\le G_\tau$ is a finite cyclic group of order at
most $3$.

(b) The map $\psi$ is injective but not surjective.

(c) The map $\psi$ is surjective but not injective and the pair
$(\psi(a),\psi(b))$ is conjugate in $G_\tau$ to $(a,b)$ or
$(a,b^{-1})$.

\item If $\sigma, \tau  \in U_m$ and $|\sigma |<|\tau|$ then for every
  homomorphism $\psi:G_\sigma \to G_\tau $ the image $\psi(G_\sigma)$ is a
  finite cyclic group of order at most $3$.

\item If  $\sigma, \tau\in U_m$ and $|\sigma|=|\tau|$ then for any
  homomorphism $\psi:G_\sigma\to G_\tau$ either the image $\psi(G_\sigma)$ is a
  finite cyclic group of order at most $3$ or $\psi$ is an
  isomorphism.

\item For every $\tau\in Q_m$ the group $G_\tau$ is a complete, Hopfian  and co-Hopfian.
\end{enumerate}

\end{theor}

  Since  $G$ is a complete, Hopfian, co-Hopfian group, if $\psi$
is any endomorphism of $G$ then  $\psi$ is injective if and only
if  $\psi$ is surjective if and only if  $\psi$ is an inner
automorphism.  It seems likely that ``endomorphism rigidity''
is another general aspect  of ``randomness'':
A random structure should not have any  endomorphisms
except  those absolutely required by the  nature of the structure.

 The Hopficity and co-Hopficity of $G_\tau$ in Theorem~\ref{B}
are the analog of  deep  results of Sela about torsion-free hyperbolic
groups~\cite{Se1,Se2}.   But the proofs here are
much simpler precisely because of the limited torsion  and do not require  Rips' machinery for analyzing group
actions on $\mathbb R$-trees.

  In considering quotients of the modular group we consider groups 
which are 
obviously  presented
as such quotients, that is, presentations of the form
\[ G = \langle a, b\  |\  a^2=b^3=1, \ r_1=1,\dots , r_m=1 \rangle  .\]
Schupp \cite{S} proved that the triviality problem restricted to
such presentations remains undecidable.  The isomorphism problem
for such presentations is thus certainly undecidable.  Indeed, the
proof provided shows that the isomorphism problem restricted
to a  fixed class $\mathcal P_m$ defined immediately below is 
undecidable for all $m \ge 15$. (This uses the fact that there is
a $2$-generator, $11$-relator group with unsolvable word problem.)
Nonetheless, rigidity shows that the isomorphism problem is
generically easy (see \cite{KMSS} for the definitions
   of generic-case complexity):

\begin{cor}\label{generic complexity}
   Let $m\ge 1$ be an arbitrary integer. Let $\mathcal P_m$ be the
   class of presentations of the form
   \[
\langle a,b\ |\  a^2=b^3=1, r_1=1,\dots, r_m=1\rangle,
   \]
   where $r_i$ are cyclically reduced words in $M$.

   Then the isomorphism problem for groups defined by
   presentations from the class $\mathcal P_m$ is strongly
   generically quartic time.
\end{cor}

\begin{proof}
We will describe a partial algorithm that solves the isomorphism
problem for $\mathcal P_m$ strongly generically in quartic time.

   The set of pairs of tuples $(\sigma, \tau)\in \mathcal P_m\times \mathcal P_m$ such that one of
$R(\sigma), R( \tau)$ satisfies the  $Q_m$ condition and the other
satisfies the standard $C(1/8)$ small cancellation condition is
exponentially generic. By Theorem~\ref{A} we can verify if this is
indeed the case in quartic time.

If it is not the case, the algorithm does not return any answer.
If the condition is satisfied then by Theorem~\ref{A}, we know
that $G_{\tau}$ is isomorphic to $G_{\sigma}$ if and only if
 $R(\tau)=R(\sigma)$ or  $R(\tau)=\eta (R(\sigma))$.  We can verify
if one of these equalities holds in cubic time. If it does the
groups are isomorphic and if not then the groups are not
isomorphic.
\end{proof}

  Once   one has rigidity one can compute  the exact  asymptotics of the
number of \emph{isomorphism  types}
of groups given by   relevant  presentations.  Note that the statement of Theorem~\ref{C}
below does not involve the notion of genericity in any way.

\begin{theor}\label{C}[Counting Isomorphism Types]
  Let $m\ge 1$. Let $I_m(n)$ denote the number of isomorphism types of
  groups given by presentations of the form
\[
M/\langle\langle r_1,\dots,r_m\rangle\rangle,
\]
where $r_1,\dots, r_m$ are cyclically reduced words of length $n$ in
$M$.

Then for even $n\to\infty$
\[
I_m(n)\sim \frac{(2^{\frac{n}{2}+1})^{m}}{2\ m!(2n)^m},
\]
that is
\[
\lim_{k\to\infty}\frac{2 I_{m}(2k)\ m!(4k)^m}{(2^{k+1})^{m}}=1.
\]
\end{theor}

  The  theory of Kolmogorov complexity is a general theory of
``descriptive complexity'' and the first basic result is that
a long random word over a finite alphabet is essentially its own shortest description.
One might summarize this result by saying that random words are ``essentially incompressible''.
Rigidity is inherent in this situation since two words are equal only if
they are identical.  In  other  situations  where there is   rigidity
for algebraic structures,  one can also investigate the appropriate descriptive complexity.
In the case of groups, the idea of  the $T$-invariant  was introduced by Delzant~\cite{De}.
We need here  a slight variation which we call the  \emph{absolute} $T$-invariant.

\begin{defn}\label{defn:T1}[Absolute $T$-invariant]
  Let
\[
\Pi=\langle a_1,\dots, a_s| w_1,\dots, w_t\rangle
\]
be a finite group presentation.  We define
$\ell_1(\Pi):=\sum_{i=1}^t |w_i|$.

For a finitely presentable group $G$ let $T_1(G)$ be the minimum
of $\ell_1(\Pi)$ taken over all finite presentations $\Pi$ of $G$.
The number $T_1(G)$ is called the \emph{absolute $T$-invariant} or
the \emph{ descriptive complexity} of $G$.
\end{defn}

The definition of $T_1(G)$ differs slightly from
\emph{Delzant's $T$-invariant $T(G)$} of a finitely presentable
group $G$ where  the ``length''
 being minimized is
$\ell(\Pi)=\sum_{i=1}^t \max\{0, |w_i|-2\}$. It turns out that
$T(G)$ is better for certain   topological arguments. In
particular, Delzant proved that $T(G_1\ast G_2)=T(G_1)+T(G_2)$.
The $T$-invariant plays an important role in Delzant and
Potyagailo's proof of the strong accessibility (or "hierarchical
decomposition") theorem for finitely presented groups~\cite{DP}.
Both $T(G)$ and $T_1(G)$ are related to the notion of Matveev
complexity for 3-manifolds, and this connection is explored in a
recent paper of Pervova and Petronio~\cite{PP}.

\begin{theor}\label{D}[Essential Incompressibility]
Let $m\ge 1$ be a fixed integer.  For any $0<\epsilon<1$ there is an
  integer $n_0>0$ and a constant $L=L(m,\epsilon)>0$ with the
  following property.

  Let $J$ be the set of all tuples $\tau\in \mathcal T_m$ such that
\[
T_1(G_\tau)\log_2 T_1(G_\tau) \ge L |\tau|.
\]
Then for any $n\ge n_0$

\[
\frac{\gamma(n,J)}{\gamma(n,\mathcal T_m)}\ge 1-\epsilon.
\]
\end{theor}

Informally, Theorem~\ref{D} says that for an ``almost generic''
$m$-relator quotient $G_\tau$ of $M$ with relators of equal length
the function $T_{1}( G_\tau)$ is bounded from below by an ``almost
linear''  function in terms of the length of the given
presentation of $G_\tau$. Thus ``almost generic'' $m$-relator
quotients of $M$ are  \emph{essentially  incompressible} in the
sense that the given presentation is almost the shortest possible
description of the group. While such a conclusion may  not be
unexpected, it is surprising that one is actually able to prove
this result.  Note that $T_{1}(G) = 0$ if and only if $G$ is a
free group and that the only free quotient of the modular group is
the trivial group.   Since the triviality problem is undecidable
for quotients of the modular group \cite{S}, we cannot in general
even decide if $T_{1} = 0$   Theorem~\ref{D} is  a generalization
to the present situation of a similar result for random
one-relator groups obtained by the authors in \cite{KSn}. As to be
expected, basic results about  Kolmogorov complexity are crucial
to the proof of Theorem~\ref{D}.

Our previous  results on  rigidity and related topics were
restricted to one-relator groups ~\cite{KS,KSn,KSS} because the
arguments relied on a classic result of Magnus~\cite{Magnus}: If
two elements $r$ and $s$ in a free group $F$ have the same normal
closures then $r$ is conjugate to $s^{\pm 1}$ in $F$. Similar
statements are generally false for subsets of free groups with
more than one element. In this paper we overcome that  difficulty
in studying generic quotients of the modular group.  First note
that the set of infinite quotients of the modular group  is  in
some sense ``very ample''.
 Miller and Schupp ~\cite{MS}  showed that every
countable group  can be embedded in a complete, Hopfian quotient
of the modular group and  the embedding preserves the property of
being finitely presented.  Schupp ~\cite{S} later proved that
every countable group can be embedded in a simple group which is
quotient of the modular group. The  presence of torsion of a very
restricted nature greatly limits homomorphisms between  quotients
of $M$ and, like the papers cited above, we exploit that idea
here. Given the extra control on homomorphisms, it turns out to be
possible to replace the result of Magnus mentioned above by its
general analog in small cancellation theory.
Greendlinger~\cite{Green} proved the following theorem: If $S$ and
$R$ are symmetrized $C'(1/6)$-subsets of a free group  with
$\langle\langle S\rangle\rangle=\langle\langle R\rangle\rangle$
then $S=R$. The proof of the similar  result is essentially unchanged
for subsets of $M$ satisfying a suitable small cancellation
condition.

We adopt the ``Arzhantseva-Ol'shanskii method'' to study quotients
of the modular group.  First,  we shall use a version of the very strong
``nonreadability'' small cancellation hypothesis introduced by
Arzhantseva and Ol'shanskii ~\cite{AO}. Verifying that $m$-relator
presentations over the modular group which satisfy
 such a condition is a  generic set is simpler than for
quotients of free groups.  Second,  we represent subgroups of $M$
by labelled graphs as usual. In our case, these are
$A$-\emph{graphs}, that is graphs where edges are labelled by the
elements of $A = \{ a,b,b^{-1} \} \subseteq M$.  To take advantage
of the strong small cancellation condition one needs to perform
certain \emph{Arzhantseva- Ol'shanskii moves} which  preserve the
subgroup of $G$  represented by an $A$-graph $\Gamma$. A crucial
observation established in the present paper is that performing
such moves on \emph{separating} arcs results in \emph{Torelli
equivalence} at the level of generating tuples.

There are several places in the proofs of the main results of this
paper where we substantially use the fact that we are working
specifically with the quotients of $M$. First, it is important for
our arguments to know that $M$ is a free product of \emph{finite}
cyclic groups. We use small cancellation considerations to
conclude that for a generic quotient $G=M/N$ of $M$ every element
of order $2$ is conjugate to $a$ in $G$ and every element of order
$3$ is conjugate to $b^{-1}$. Thus we know that if $\psi:M\to G$
is a homomorphism with $\psi(a)\ne 1$ and $\psi(b)\ne 1$ then
$(\psi(a),\psi(b))=(u_1au_1^{-1},u_2b^{\pm 1}u_2^{-1})$. Next we
need to conclude that if $\psi$ as above is not injective then in
fact the 2-tuple $(\psi(a),\psi(b))$ is conjugate in $G$ to
$(a,b^{\pm 1})$. At this point we use $AO$-moves and genericity
assumptions to prove that $(\psi(a),\psi(b))$ is \emph{Torelli
equivalent} and hence conjugate to $(a,b^{\pm 1})$. This is the
place in the proof where it is critically important for the
argument that the number of generators be equal two $2$. Indeed,
for 2-tuples of elements Torelli equivalence is the same as
conjugacy. However, for $k$-tuples with $k\ge 3$ Torelli
equivalence is not the same thing as conjugacy and our arguments
do not apply in that case. The reason for this is a well-known
fact that for $k\ge 3$ the \emph{Torelli subgroup} of $Aut(F_k)$,
consisting of all automorphisms of $F_k$ that induced the identity
map in the abelianization of $F_k$, is strictly bigger than the
group of inner automorphisms of $F_k$.

After establishing that $(\psi(a),\psi(b))$ is conjugate to
$(a,b^{\pm 1})$ in $G$, we apply Greendlinger's Theorem that
symmetrized small cancellation sets with the same normal closures
are equal, to conclude the proof of Theorem~\ref{A}.

We strongly believe that the results of this paper should also
hold for generic quotients of free groups of arbitrary finite rank
$k\ge 2$, and, more generally, of free products of $k\ge 2$ cyclic
groups. Carrying out a proof in that context requires establishing
some version of the following ``Stability Conjecture''. Our computer 
experiments so far tend to support the conjecture.

\begin{conj}[The Stability  Conjecture]
  Fix $k\ge 2$ and $m \ge 1$ and let $F=F(a_1,\dots, a_k)$. Then there
  exists an algorithmically recognizable generic class $\mathcal{Y}$
  of $m$-tuples of elements of $F$ with the following property. If
  $\sigma,\tau \in \mathcal{Y}$ and $\alpha\in Aut(F)$ are such that
  $R(\sigma)$ and $R(\alpha(\tau))$ have the same normal closure in
  $F$ then $R(\sigma) = R(\alpha(\tau))$.
\end{conj}

\section{Arzhantseva-Ol'shanskii moves on graphs and Torelli equivalence}

      Since we represent subgroups by labelled graphs we completely list our conventions.
  We follow the same conventions about graphs as Serre.

\begin{conv}

   A \emph{graph} is   a tuple $\Gamma=(V,E, o,t,{}^{-1})$
  where $V=V(\Gamma)$ is the \emph{vertex set} of $\Gamma$, $E=E(\Gamma)$
  is the \emph{edge set} of $\Gamma$ and $o:E\to V$, $t:E\to V$ and
  ${}^{-1}:E\to E$ are the \emph{origin, terminus} and \emph{inverse} maps.
  We require that ${}^{-1}:E\to E$ is  an involution
   with  $e^{-1}\ne e$ and  $t(e)=o(e^{-1})$.

  An \emph{orientation} on $\Gamma$ is a partition $E=E^+\sqcup E^-$
  such that  $e\in E^+$ if and only if
  $e^{-1}\in E^-$.

  An \emph{arc}  is a simple edge-path where all vertices
  have degree 2 in $\Gamma$, except possibly for the initial and the
  terminal vertices.

\end{conv}

   We  discuss a version of Arzhantseva-Ol'shanskii
moves in the general setting of abstract graphs. The key
observation is that performing such a move on a separating arc of
a graph naturally corresponds to the Torelli equivalence at the
level of the free bases of the fundamental groups.

\begin{defn}[Abstract $AO$-move]\label{defn:abstrAO}
  Let $\Gamma$ be a connected graph. Let $p_1\gamma p_2$ be a path in
  $\Gamma$ such that $\gamma$ is a non-loop arc of $\Gamma$ and the
  paths $p_1$, $p_2$ do not pass through $\gamma$ or $\gamma^{-1}$.

  Modify $\Gamma$ by first attaching to $\Gamma$ a new arc $\gamma'$
  (possibly consisting of several edges) from $o(p_1)$ to $t(p_2)$ and
  then removing the arc $\gamma$. The resulting graph $\Gamma'$ is
  said to be obtained from $\Gamma$ by \emph{a  move of type
    $AO$}.  We define the \emph{$AO$-map $P:\Gamma\to\Gamma'$
    associated to this  $AO$-move} as follows.  We set $P$ to
  be the identity map on all edges and vertices of $\Gamma$ that are
  not changed by the $AO$-move and thus are common for $\Gamma$ and
  $\Gamma'$. This includes the endpoints of $\gamma$ and $\gamma'$. We
  define $P$ on $\gamma$ to "push" $\gamma$ to the path
  $p_1^{-1}\gamma'p_2^{-1}$. The map $P:\Gamma\to\Gamma'$ is a
  homotopy equivalence. Indeed, undoing the $AO$ move from $\Gamma$ to
  $\Gamma'$ is an $AO$-move from $\Gamma'$ to $\Gamma$ consisting in
  removing $\gamma'$ and adding back $\gamma$.  The map
  $P':\Gamma'\to\Gamma$ defined similarly to $P$ is easily seen to be
  a homotopy inverse of $P$. Thus $P:\Gamma\to\Gamma'$ is a homotopy
  equivalence.
\end{defn}

\begin{defn}[Elementary Torelli moves]
  Let $\tau=(g_1,\dots, g_n)$ be an $n$-tuples of elements of a group $G$.
The   following transformations of $\tau$ will be called the
  \emph{elementary Torelli moves}:
\begin{enumerate}
\item For some subset $S\subseteq \{1,2,\dots, n\}$ and some $h\in
  \langle \{g_j| j\in S\}\rangle\le G$ for each $i\in S$ replace $g_i$
  by $hg_ih^{-1}$.
\item For some subset $S\subseteq \{1,2,\dots, n\}$ and some $g\in
  \langle \{g_j| j\not\in S\}\rangle\le G$ for each $i\in S$ replace
  $g_i$ by $gg_ig^{-1}$.
\item Conjugate the entire tuple $\tau$ by some element $g\in G$.
\end{enumerate}
\end{defn}

\begin{defn}[Torelli Equivalence]
  We say that two $n$-tuples $\tau$ and $\tau'$ of elements of $G$ are
  \emph{Torelli-equivalent} if there
  exists a finite chain of Torelli moves taking $\tau$ to $\tau'$.
\end{defn}

Note that if $n=2$ we have an ordered pair of elements and  then any two Torelli-equivalent pairs
generate conjugate subgroups of $G$.

\begin{notation}
  Let $\Gamma$ be a finite connected graph free of rank $n\ge 1$. Let
  $T$ be a maximal subtree of $\Gamma$.

  If $x_0\in V\Gamma$ is a base-vertex, then $T$ defines a free basis
  of $\pi_1(\Gamma, x_0)$ as follows. For each edge $e\in
  E^+(\Gamma-T)$ define the path $c(e)$ as $c(e):=[x_0,o(e)]_T e
  [t(e),x_0]$, where $[x,y]_T$ stands for the unique reduced edge-path
  from $x$ to $y$ in $T$.

  Choose an orientation $E\Gamma=E^+\Gamma\cup E^{-}\Gamma$ and choose
  an ordering $e_1,\dots, e_n$ of all the elements of $E^+(\Gamma-T)$.
  Then the $n$-tuple $(c(e_1),\dots, c(e_n))$ is a free basis of
  $\pi_1(\Gamma, x_0)$.
  We will denote the tuple $(c(e_1),\dots, c(e_n))$
  by $S_{T,x_0}$. While $S_{T,x_0}$ does depend on the choice of an
  orientation on $E(\Gamma-T)$ and on the choice of an ordering on
  $E^+(\Gamma-T)$, these choices will usually be fixed and explicit
  references to them will be suppressed.
\end{notation}

\begin{conv}[$AO$ move on a separating arc]\label{conv:ord}
Let $\Gamma'$ be obtained from a finite connected graph $\Gamma$
by an abstract $AO$-move  via removing an arc $\gamma$ and adding
an arc $\gamma'$,
 as in Definition~\ref{defn:abstrAO}.
Suppose that $\gamma$ is a separating arc of $\Gamma$. Let $T$ be a
maximal subtree in $\Gamma$, so that $T$ contains $\gamma$.   Suppose
that we fix an orientation on $\Gamma$ and an ordering on $E^+(\Gamma-T)$.

Put $T':=T-\{\gamma\}
\cup \{\gamma'\}$. Thus $T'$ is a maximal subtree in $\Gamma'$.

Note that $E(\Gamma-T)=E(\Gamma'-T')$. We choose an orientation on
$\Gamma'$ so that it agrees with $\Gamma$ on $\Gamma'-T'$, so that
$E^+(\Gamma-T)=E^+(\Gamma'-T')$. We also order $E^+(\Gamma'-T')$
exactly as $E^+(\Gamma-T)$ was ordered.

Note also that $x_0=o(\gamma)=P(x_0)$ is a vertex of both $\Gamma$ and $\Gamma'$

The tuples $S_{T,x_0}$ and $S_{T',x_0}$ of elements of
$\pi_1(\Gamma,x_0)$ and of $\pi_1(\Gamma',x_0)$ are defined according
to the above conventions regarding the orientations and the orderings
of positive edges outside of the specified maximal trees.
\end{conv}

In general, $A0$-moves on \emph{nonseparating arcs} result in Nielsen
equivalence at the level
of the tuples $S_{T,x_0}$ (see~\cite{A4,KS}). We observe here that
if an $AO$-move is performed on a \emph{separating} arc $\gamma$,
this move results  in Torelli
equivalence:

\begin{prop}\label{tor}
  Let $\Gamma$ be a finite connected graph with the fundamental group
  free of rank $n\ge 1$ and without degree-one vertices.  Let
  $\Gamma'$ be obtained from $\Gamma$ by an abstract $AO$-move
  removing an arc $\gamma$ and adding an arc $\gamma'$.  Let
  $P:\Gamma\to\Gamma'$ be the $AO$-map corresponding to this
  $AO$-move.  Let $x_0=o(\gamma)$, so that $P(x_0)=x_0$.

  Suppose that $\gamma$ is a separating arc of $\Gamma$.

  Let $T$ be a maximal subtree in $\Gamma$ (and hence $T$ contains
  $\gamma$).  Put $T':=T-\{\gamma\} \cup \{\gamma'\}$. Thus $T'$ is a
  maximal subtree in $\Gamma'$. Let the orientations and the orderings
  of positive edges outside of maximal trees are chosen as in
  Convention~\ref{conv:ord}.

  Let  $P_\#: \pi_1(\Gamma, x_0)\to \pi_1(\Gamma', x_0)$
be the homomorphism of fundamental groups induced.
by the $AO$-map $P$ above.
  Then the tuples $P_\#(S_{T,x_0})$ and $S_{T',x_0}$ are
  Torelli-equivalent in $\pi_1(\Gamma',x_0)$.
\end{prop}

\begin{proof}

  Since $\gamma$ is a separating arc and is thus  not a loop, the graph
  $\Gamma-\{\gamma\}$ consists of two connected components: $\Gamma_1$
  containing $x_0=o(\gamma)$ and $\Gamma_2$ containing $t(\gamma)$.
  The set $E^+(\Gamma-T)$ is partitioned as: $\{e_1,\dots, e_k,
  f_1,\dots, f_{n-k}\}$ where $e_i\in \Gamma_1$ and $f_j\in \Gamma_2$.
  Moreover, $T_q=T\cap \Gamma_q$ is a maximal tree in $\Gamma_q$ for
  $q=1,2$.

     For future reference we need to explicitly write down the
 $n$-tuple $S_{T,x_0}$ corresponding to $T$.
   Let $y:=[o(\gamma), o(p_1)]_T$ and
  let $z:=[t(\gamma),t(p_2)]_T$.  For each $e_i\in E^+(\Gamma_1-T)$
  let $y_i:=[o(\gamma),o(e_i)]_T$ and let $y_i':=[t(e_i), o(\gamma)]_T$.
  Similarly for each $f_j$ let $z_j$ and $z_j'$ be the paths
  $[t(\gamma), o(f_j)]_T$ and of $[t(f_j),t(\gamma)]_T$ accordingly.
  Note that $[t(\gamma'),o(f_j)]_T$ is homotopic relative endpoints in
  $T$ to $z^{-1}z_j$ and that $[t(f_j),t(\gamma')]_T$ is homotopic
  relative endpoints in $T$ to $z_j'z$.

  Recall that $x_0=o(\gamma)$ is the base-vertex of $\Gamma$
 and note that $x_0=P(x_0)$ is still the base-point of
$\Gamma'$.

  Then by definition the $n$-tuple $S_{T,x_0}$ corresponding to $T$
  has the form:
\[
S_{T,x_0}=(r_1,\dots, r_k, s_1, \dots, s_{n-k}),
\]
where $r_i= y_i e_i y_i'$ and $s_j=\gamma z_jf_jz_j'\gamma^{-1}$.

Let $\Gamma'$ now be obtained from $\Gamma$ by removing $\gamma$ and
adding an arc $\gamma'$ from $o(p_1)$ to $t(p_2)$.

By definition of $P$ we have $P_\#(r_i)=r_i$ for $i=1,\dots, k$. Also
for $j=1,\dots, n-k$ we have
\[
P_\#(s_j)=p_1^{-1}\gamma'p_2^{-1}z_jf_jz_j^{-1}p_2(\gamma')^{-1}p_1.
\]

Recall that $T':=T-\{\gamma\} \cup \{\gamma'\}$ is a maximal subtree
in $\Gamma'$ and We wish to compute explicitly the tuple $S_{T',x_0}$.

Note that $E^+(\Gamma-T)=E^+(\Gamma'-T')=\{e_1,\dots, e_k, f_1,\dots,
f_{n-k}\}$.

Clearly, the elements of $S_{T',x_0}$ corresponding to $e_i$ remain
the same for $\Gamma'$ as they were for $\Gamma$, that is $r_i= y_i
e_i y_i'$. The elements $s_j$ corresponding to the $f_j$ will change to

\begin{gather*}
  s_j'=y\gamma' [t(\gamma'), o(f_j)]_T f_j[t(f_j),t(\gamma')]_T (\gamma')^{-1}y^{-1}\\
  =y\gamma' z^{-1}z_jf_jz_j'z (\gamma')^{-1}y^{-1}.
\end{gather*}

Since $y p_1$ is a loop at $x_0$ in $\Gamma_1$ and $y\gamma'
p_2^{-1}z(\gamma')^{-1}y^{-1}$ is a loop at $x_0$ in $ y\cup
\gamma'\cup \Gamma_2$, we conclude that $y p_1=W_1(r_1,\dots, r_k)$ in
$\pi_1(\Gamma,x_0)$ and
$y\gamma'p_2^{-1}z(\gamma')^{-1}y^{-1}=W_2(s_1',\dots, s'_{n-k})$ in
$\pi_1(\Gamma,x_0)$ for some words $W_1,W_2$.

Recall also that $P(\gamma)=p_1^{-1}\gamma'p_2^{-1}$.  Therefore in
$\pi_1(\Gamma',x_0)$ we have

\begin{align*}
  &P_\#(s_j)=p_1^{-1}\gamma'p_2^{-1} z_j f_j z_j' p_2(\gamma')^{-1}
  p_1\\
  &=(p_1^{-1}y^{-1}y\gamma' p_2^{-1}z(\gamma'^{-1})y^{-1}) (y\gamma'
  z^{-1}z_j f_j z_j'
  z (\gamma')^{-1}y^{-1}) (y\gamma'z^{-1} p_2(\gamma')^{-1}y^{-1}yp_1 )=\\
  &=W_1 W_2 s_j' W_2^{-1} W_1^{-1}.
\end{align*}

Recall that $W_1=W_1(r_1,\dots, r_k)$ and $W_2=W_2(s_1',\dots,
s_{n-k}')$. Hence the $n$-tuples
\[
S_{T',x_0}=(r_1, \dots, r_k, s_1', \dots, s_{n-k}')
\]
and
\[
P_\#(S_{T,x_0})=(r_1, \dots, r_k, W_1 W_2 s_1' W_2^{-1} W_1^{-1},
\dots, W_1 W_2 s_{n-k}' W_2^{-1} W_1^{-1})
\]
are Torelli-equivalent in $\pi_1(\Gamma',x_0)$, as claimed.
\end{proof}

\begin{conv}
An \emph{$A$-graph} is a graph $\Gamma$ where every edge $e$
  is labelled by an element $\mu(e)\in A$ such that
  $\mu(e^{-1})=(\mu(e))^{-1}$. Thus
  $\mu(e)=a$ iff  $\mu(e^{-1})=a$ and $\mu(e)=b$ iff
  $\mu(e^{-1})=b^{-1}$.  For an edge-path $p = e_1 ... e_l$ in $\Gamma$ the
  \emph{label} $\mu(p) = \mu(e_1) ... \mu(e_l)$ is a word in the alphabet $A$.

Let $G=M/N$ be a fixed quotient of the modular group.
  If $\Gamma$ is an $A$-graph with a base-vertex $x_0$, there is
  a natural \emph{labelling homomorphism} $\phi: \pi_1(\Gamma,x_0)\to
  G$ that sends the homotopy class of a  closed edge-path at $x_0$ to
  the element of $G$ represented by the label of that path. We say
  that $H=\phi(\pi_1(\Gamma,x_0))\le G$ is the \emph{subgroup of $G$
    represented by $(\Gamma,x_0)$}. If $\Gamma$ is connected
  then the conjugacy class of $H$ does not depend on the choice of the
  base-vertex.
\end{conv}

\begin{defn}[Arzhantseva-Ol'shanskii move: move $AO$ on $A$-graphs]\label{defn:AO}
  Let $N\triangleleft M$ and $G=M/N$ be a fixed quotient of $M$.
  Suppose $\Gamma$ is a connected $A$-graph. Let $p_1\gamma p_2$ be a
  path in $\Gamma$ such that $\gamma$ is a non-loop arc of $\Gamma$
  and the paths $p_1$, $p_2$ do not pass through $\gamma$ or
  $\gamma^{-1}$. Let $u_1,u_2$ be the labels of $p_1,p_2$ and let $u$
  be the label of $\gamma$. Suppose $v$ is a reduced word in $A$ such that
  $u_1uu_1=v$ in $G$.

  Modify $\Gamma$ by first attaching to $\Gamma$ a new arc $\gamma'$
  labelled $v$ from $o(p_1)$ to $t(p_2)$ and then removing the arc
  $\gamma$. The resulting $A$-graph $\Gamma'$ is said to be obtained from
  $\Gamma$ by \emph{an $AO$-move}.
\end{defn}

Proposition~\ref{tor} immediately implies:

\begin{cor}\label{tech}
  Let $N\triangleleft M$ and $G=M/N$ be a fixed quotient of $M$.
  Let $\Gamma$ be a finite connected $A$-graph with the fundamental
  group free of rank $n\ge 1$.  Suppose that an $AO$-move applies to
  $\Gamma$ and let $p_1$, $\gamma$, $p_2$, $\gamma'$, $\Gamma'$ be as
  in Definition~\ref{defn:AO}.

  Suppose also that $\gamma$ is a separating arc of $\Gamma$.

  Let $T$ be a maximal subtree in $\Gamma$ (and hence $T$ contains
  $\gamma$).  Put $T':=T-\{\gamma\} \cup \{\gamma'\}$. Thus $T'$ is a
  maximal subtree in $\Gamma'$.

  Let $x_0\in V\Gamma$ and $x_0'\in V\Gamma'$ be base-vertices. Let
  $\phi: \pi_1(\Gamma, x_0)\to G$, $\phi': \pi_1(\Gamma, x_0')\to G$
  be the labelling homomorphisms.

  Then the tuples $\phi(S_{T,x_0})$ and $\phi'(S_{T',x_0'})$ are
  Torelli-equivalent in $G$.

  In particular, if $n=2$, then the 2-tuples $\phi(S_{T,x_0})$ and
  $\phi'(S_{T',x_0'})$ are conjugate in $G$.

\end{cor}

\section{The Generic Nonreadability Condition}

   Recall that we are using  the group alphabet $A=\{a,b, b^{-1}\}\subseteq M$.
and that a   word $w$ in the alphabet $A$ is \emph{reduced} if it
does not contain any subwords of the form $aa,bb^{-1}, b^{-1}b,bb,
b^{-1}b^{-1}$.
Note that if $w\in \mathcal C$ is a cyclically reduced word with
$|w|>1$ then either $w$ begins with $a$ and ends with $b^{\pm 1}$ or
$w$ begins with $b^{\pm 1}$ and ends with $a$  and so
$w$ has even length  in both cases.

The following lemma is therefore straightforward.

\begin{lem}\label{est} The following hold:
\begin{enumerate}
\item If $n\ge 2$ is an even integer then
\[
\gamma(n,\mathcal C)=\gamma(n+1,\mathcal C)=2\cdot 2^{n/2}.
\]
\item There exist constants $c_1,c_2>0$ such that for every $n\ge 1$
\[
c_1 2^{n/2}\le \rho(n,\mathcal C)\le c_2 2^{n/2}
\]
\end{enumerate}
\end{lem}

  The following statement is a
straightforward corollary of the definitions of genericity and of
Lemma~\ref{est}:

\begin{prop}\label{prob:basic} Let $m\ge 1$ be an integer. Then the following hold:
\begin{enumerate}
\item A subset $S\subseteq \mathcal C^m$ is exponentially negligible
  in $\mathcal C^m$ if and only if
\[
\lim_{n\to\infty} \frac{\rho(n,S)}{2^{mn/2}}=0
\]
and the convergence is exponentially fast.
\item Let $S\subseteq \mathcal T_m$. Then $S$ is exponentially
  negligible in $\mathcal C^m$ if and only if $S$ is exponentially
  negligible in $\mathcal T_m$.

\item A subset $S\subseteq \mathcal T_m$ is exponentially
  negligible in $\mathcal T_m$ if and only if for even $n\to\infty$

\[
\frac{\gamma(n,S)}{2^{mn/2}}\to 0
\]
and the convergence is exponentially fast.
\end{enumerate}
\end{prop}

   The idea of an  Arzhantseva-Ol'shanskii condition is that large parts of
the  relators which  one wants to consider are \emph{not} ``readable''
along certain graphs. In studying quotients of the modular group we need only
consider one type of graph.

\begin{defn} [Barbell graphs]
 Let $u$ be a word in $A^*$.  The $u$-\emph{barbell} is the graph $\Gamma$
with a loop-edge  $e$ labelled $a$, a loop-edge $f$ labelled $b$ and
 a simple arc $p$ labelled by $u$ connecting the vertex of $e$ to the
vertex of $f$. A  barbell graph $\Gamma$ is \emph{reduced} if $u$ is a cyclically
reduced word of length $2k \ge 2$ which begins with $b^{\pm1}$ and ends
with $a$.  A word $w$ is \emph{readable} in $\Gamma$ if there exists
a vertex $v \in \Gamma$ and a path $\gamma$ starting at $v$ with label $w$.
\end{defn}

   Even on a very  simple example such as  $u = ba$, one quickly
sees that the number of all  words readable on the $u$-barbell is
exponentially falling behind the number of all words.
We first need to estimate the number of all words of length $n$ readable in
a fixed reduced barbell graph.

\begin{lem}\label{lem:u}
  Let $\Gamma$ be a reduced $u$-barbell graph where $|u|=2k$ and $k\ge 1$.
  Then for $n\ge 1$ the number of all reduced words of length $n$
  readable in $\Gamma$ is at most
\[
c_3 2^k 2^{n/(4k+2)},
\]
where $c_3>1$ is a constant independent of $u,n,k$.
\end{lem}
\begin{proof}

  The graph $\Gamma$ has $(2k + 1)$ vertices and for each vertex $v$
on can begin reading in either of two directions.  Let $w$ be any word of
  length $n$ readable in $\Gamma$. Then $w$ can be written as
\[
w=w_1w'w_2
\]
where each $|w_i| < (4k + 1)$ and transverses the $b$-loop at most once,
and where $w'$ has the form
\[
w'=(aub^{\epsilon_1}u^{-1})\dots (aub^{\epsilon_{t}}u^{-1})
\]
with $\epsilon_j\in \{1, -1\}$.  The number of possibilities for such a $w'$
is  $2^{t} \le  2^{n/(4k+2)} $ and so the number of all
possibilities for $w$ is at most
\[
 2 \cdot 2 \cdot (2k+ 1)  2^{n/(4k+2)} \le c_3 2^k  2^{n/(4k+2)}
\]
for  $c_3 2^k \ge 8k + 4$  so $c_3$ is a  constant independent of $u,n,k$.

\end{proof}

\begin{lem}\label{lem:c4}
  Let $0< \theta < \frac{1}{30}$.  Then the number of reduced words of
  length $n\ge 1$ that are readable as labels of some paths in
  reduced $u$-barbell graphs for some $u$ with $|u|\le \theta n$ is at most
\[
c_4 2^{n/5}
\]
where $c_4>0$ is some constant independent of $n$ and $\theta$.
\end{lem}
\begin{proof}
Let $0<2k\le \theta n$.  The number of reduced words $u$ of length $2k$
which  start with $b^{\pm 1}$ and end with $a$ is equal to $2^k$.

Therefore by Lemma~\ref{lem:u} the number of all reduced words of length $n\ge 1$ that are readable in
$u$-barbell graphs with $|u|\le \theta n$ is at most
\begin{gather*}
  c_3\sum_{2\le 2k\le\theta n} 2^k 2^k  2^{n/(4k+2)}=c_3\sum_{2\le
    2k\le\theta n}
  2^{2k} 2^{n/(4k+2)}\le \\
  c_3\sum_{2\le 2k\le\theta n}2^{2\theta n/2} 2^{n/(4k+2)}\le
  c_3\sum_{2\le 2k\le\theta n}2^{\theta n} 2^{n/6}\le\\
  \le c_3\frac{\theta n}{2} 2^{\theta n + n/6} \le c_4 2^{n/5}
\end{gather*}
where $c_4>0$ is some constant independent of $\theta, n$ and where the
last inequality holds since by the choice of $\theta$ we have
\[
\frac{n}{6}+\theta n <\frac{n}{5}.
\]
\end{proof}

\begin{defn}
Let $0<\theta \le  \frac{1}{40}$. A cyclically reduced word $w$ is said to
be \emph{$\theta$-readable} if there is a subword $v$ of some cyclic
permutation of $w$ or $w^{-1}$ with $|v|\ge |w|/2$ such that $v$ is
readable in some reduced $u$-barbell graph with $|u|\le \theta |v|$.
\end{defn}

\begin{defn}[Genericity conditions]\label{defn:Q}
  Let $0<\lambda\le \frac{1}{120}$. Let $m\ge 1$ be an integer. We say
  that a tuple $(r_1,\dots, r_m)$ of cyclically reduced words in
  $M$ satisfies the \emph{$Q_m(\lambda)$-condition} if the
  following hold:
\begin{enumerate}
\item For each $i=1,\dots, m$ the word $r_i$ is not $\frac{1}{40}$-readable.
\item The symmetrized closure of $\{r_1,\dots, r_m\}$ satisfies the
  $C'(\lambda)$-small cancellation condition in $M$.
\item For each $i=1,\dots, m$ the word $r_i$ is not a proper power in
  $M$ (in particular, this means that $r_i$ is different from
  every cyclic permutation of $r_i$).
\item If $i\ne j$ then $r_i$ is not a cyclic permutation of $r_j^{\pm
    1}$.
\item If $1\le i\le m$ then no subword $z$ of any cyclic
permutation of  $\eta(r_i)$  with  $ |z| > |\eta(r_i)|/3 $ occurs
as a subword of any  cyclic permutation of any $r_j^{\pm 1}$ for
$j=1,\dots, m$.

\item For every $1\le i\le m$ the word $r_i$ is not a cyclic
  permutation of $r_i^{-1}$.

\item For each $i=1,\dots, m$ we have $|r_i|\ge 1200$.
\end{enumerate}

    The small cancellation condition of item $(2)$ implies that
 $(4)$ and  $(6)$ hold but we have listed them for convenience.
We say that $(r_1,\dots, r_m)\in \mathcal C^m$ satisfies condition
$U_m(\lambda)$ if $(r_1,\dots, r_m)\in Q_m(\lambda)$ and, in addition,
$|r_1|=\dots=|r_m|$. Thus $U_m(\lambda)=Q_m(\lambda)\cap\mathcal T_m$.
\end{defn}

\begin{rem}
  Note that conditions $Q_m(\lambda)$ and $U_m(\lambda)$ are invariant
  under reordering the tuple, inverting member of the tuple, taking a
  cyclic permutation of a member of the tuple and applying $\eta$ to
  the entire tuple.

  Recall that $R(\tau)$ is the symmetrized set generated by the elements
in the tuple $\tau$.  If $\tau=(r_1,\dots, r_m)\in U_m(\lambda)$ and
  $|r_i|=n=2k>0$ then
\[
\#R(\tau)= 2mn.
\]
   and
\[
R(\tau)\cap \eta(R(\tau))=\emptyset.
\]

\end{rem}

\begin{prop}\label{prop:generic}
  Let $0< \lambda \le \frac{1}{120}$ and let $m\ge 1$ be an integer.  Then
  the set $Q_m(\lambda)$ is exponentially generic in $\mathcal C^m$
  and the set $U_m(\lambda)$ is exponentially generic in $\mathcal
  T^m$.
\end{prop}
\begin{proof}
  We will prove that $Q_m(\lambda)$ is exponentially generic in
  $\mathcal C^m$.  The proof that $U_m(\lambda)$ is exponentially
  generic in $\mathcal T^m$ is analogous.

  It is well-known and easy to see that parts (2), (3), (4), (6) and
  (7) of Definition~\ref{defn:Q} are define exponentially generic
  subsets of $\mathcal C^m$.
   The proof that (5) defines an exponentially generic subset of
   $\mathcal C^m$ is a little more cumbersome but it is very
   similar to the proofs of Lemmas 4.6 and 4.7 in \cite{KSn},
   where the free group case was considered. We leave the details
   to the reader.

  Since the intersection of two
  exponentially generic sets is exponentially generic, it suffices to
  prove that condition (1) of Definition~\ref{defn:Q} is exponentially
  generic, that is, that the complement of condition (1) in $\mathcal
  C^m$ is exponentially negligible in $\mathcal C^m$.

  Let $M=(r_1,\dots, r_m)$ be an $m$-tuple of cyclically reduced words
  such that $|r_i|\le n$ and that part (1) of Definition~\ref{defn:Q}
  fails for $M$. Thus there is some $r_i$ such that a subword $v$ of a
  cyclic permutation of $r_i^{\pm 1}$ with $|v|\ge |r_i|/2$ has the
  property that $v$ is readable in a $u$-barbell graph with $|u|\le
  \frac{1}{40}|v|$.

  Let $c_1,c_2>0$ be the constants provided by Lemma~\ref{est}. We may
  assume that $c_2\ge 1$.

  Suppose first that $|r_i|\le 9n/10$. Then the number of
  possibilities for $M$ is at most
  \[
  c_2^m 2^{9n/10} (2^{n/2})^{m-1}.
\]
This number is exponentially smaller, as $n$ tends to infinity, than
the number $K_m(n)$ of all $M$ among all $m$-tuples of cyclically
reduced words of length at most $n$ since $K_m(n)$ satisfies
\[
K_m(n)\ge c_1^m (2^{n/2})^{m}.
\]

Suppose now that $|r_i| > 9n/10$.   Then
$r_i$ or $r_i^{-1}$ is a cyclic permutation of a word $zz'$ where
$|z| \ge |r_i|/2 > 9n/20$ and where $z$ can be read in is readable in a $u$-barbell
graph with $|u|\le \frac{1}{60}|r_i| $. Hence $|u| < n/60 < \frac{1}{30} |z|$.
By Lemma~\ref{lem:c4} the number
of such $z$ is at most $c_4 2^{n/5}$. The number of choices for the
word $z'$ of length $|z'|\le 11n/20$ is at most $c_2 2^{11n/40}$.
The number of possibilities for $zz'$ is therefore at most $c_2 c_4  2^{n/5} 2^{11n/40} \le c_5 2^{19n/40}$.
The number of cyclic permutations of any such $r_i$ and its  inverse is
at most $2n$. Thus there are at most
$2n c_5 2^{19n/40}$ possibilities for $r_i$. Since there are at most
$m$ choices for $i$, the number of possibilities for $M$ in this case
is at most

\[
2mn c_5 c_2^{m-1} 2^{19n/40} (2^{n/2})^{m-1}.
\]

Again, this number is exponentially smaller than the number $K_m(n)$
of all $m$-tuples of cyclically reduced words of length at most $n$ in
$M$. This proves that the set of all $m$-tuples $(r_1,\dots, r_m)$ of
cyclically reduced words satisfying condition $Q_m(\lambda)$ is
exponentially generic in $\mathcal C^m$.
\end{proof}

\begin{lem}\label{lem:verify}  There is a quartic-time  algorithm which,
when given an $m$-tuple $\tau = (r_1,...,r_m)$ verifies
whether or not $\tau$ satisfies the genericity condition $Q_m (\lambda)$.
\end{lem}

\begin{proof}  We first show that the condition that each $r_i$ is not $ \frac{1}{6}$-readable
can be verified in   quartic time.  If any cyclic permutation of  $r_i$ is readable in
a $u$-barbell graph $\Gamma$ we can assume that $u$ is a subword of $r^{\pm1}$.
  There are at most
\[ 2  |r_i| \  \frac{|r_i|}{60} \   \frac{|r_i|}{120}  =  |r_i|^{3} /3600 \]
such words.  We can thus construct all possible relevant $\Gamma$ in
cubic time and for each graph we can verify if  $r_i$ is readable on
it in linear time.

  We have noted that the standard small cancellation $C' (\lambda)$, item $(2)$
of the genericity condition, already implies both items $(4)$ and
$(6)$ of the condition. Let $s$ be the sum  of all the $ |r_i|$.
For a particular $ r_i$ there are at most $2 |r_i| \le 2s |r_i|/6$
subwords of length not exceeding $\lambda  |r_i|$. Whether or not
a particular word $z$ is a subword of another word $w$ can be
verified in time linear in $|z| + |w|$.  (For example, using the
Knuth-Morris-Pratt algorithm.) Verifying  $C' (\lambda)$ thus
takes at most quadratic time $c m s^2$.

  The remaining conditions can all be verified in linear time and the
lemma holds.

\end{proof}

\section{The Subgroup Theorem}

  We can now prove the
key result determining the structure of  subgroups generated
by the images of $a$ and $b$ in random quotients of $M$.
We first need to remark that while we are writing words
 as reduced words on $a,b,b^{-1}$ and
 measure lengths accordingly,  when  considering small cancellation
quotients of the modular group we must use the theory over free
products. See Lyndon-Schupp~\cite{LS} for details. In considering a van Kampen
diagram for a word equal to the identity in a small cancellation
quotient,  the theory guarantees as usual the existence of  a
region labelled by a relator $|r|$ which  has an  interior arc
$\eta$ with label $|u|$ where
 $|u| < 3 \lambda |r|$.  There is the slight technicality that
 the arc $\eta$ might begin and/or end at
 secondary vertices.  This could mean that the portion of $r$ left
on the boundary of the whole diagram could be two letters shorter
than one might think without taking this point into consideration.
By our assumption that all defining relators have length at least $100$
it follows that $ 2 < \lambda |r|$ for any defining relator and we
add extra  factor of $ \lambda $ to the lengths of interior arcs
guaranteed by the general theory.

\begin{prop}[Subgroup Theorem]\label{prop:main}

  Let $0<\lambda\le \frac{1}{120}$, let $m\ge 1$ and let
\[
G=M/N(\tau)
\]
where the tuple $ \tau = (r_1,\dots, r_m)$ of cyclically reduced words
satisfies the $Q_m(\lambda)$ condition. Then the following holds.

Suppose $g_1,g_2\in G$ are elements of orders two and three
respectively and let $H=\langle g_1,g_2\rangle\le G$. Then either
\[
H=\langle g_1\rangle \ast \langle g_2\rangle
\]
or the pair $(g_1,g_2)$ is conjugate to the pair $(a,b)$ or
$(a,b^{-1})$ in $G$.
\end{prop}

\begin{proof}

  Let $\mathcal R$ denote the symmetrized closure of $\{r_1,\dots,
  r_m\}$.  After conjugating the pair $(g_1,g_2)$ we may assume that
  $g_1=a$ and $g_2=hb^\delta h^{-1}$ for some $h\in G$ and some
  $\delta\in\{1,-1\}$.

  Among all pairs $(a,hb^\delta h^{-1})$ conjugate to $(g_1,g_2)$
  choose the pair where $|h|_G$ is the smallest possible. If $h=1$
  then $(g_1,g_2)$ is conjugate to $(a,b^\delta)$ as required. Suppose
  now that $h\ne 1$ in $G$.

  Let $u$ be a geodesic  word representing $h$.
  By minimality, $u$ is a reduced word which ends in $a$ and begins with $b^{\pm 1}$.
  Let $\Gamma$ be the barbell graph with a segment $p$ labelled by
  $u$ joining a loop-edge $e$ labelled $a$ to a loop-edge $f$ labelled
  $b$. This means that $\Gamma$ is folded and that the label of every
  reduced path in $\Gamma$ is a reduced word in $M$, provided
  that path does not contain subpaths of the form $e'e'$ where $e'$ is
  a loop-edge. Note that $\Gamma$ has $|u|+2$ non-oriented edges.

  Suppose now that $H\ne \langle g_1\rangle \ast \langle g_2\rangle$.
  Then there exists a nontrivial closed cyclically reduced path
  $\alpha$ in $\Gamma$ with label $w$ such that $w=_G 1$ and such that
  $\alpha$ does not contain subpaths of the form $e'e'$ where $e'$ is
  a loop-edge. Note that the word $w$ is reduced in $M$ by
  assumptions on $\Gamma$. Then $\alpha$ contains a subpath $\beta$
  labelled by a word $v$ such that $v$ is a subword of some $r\in
  \mathcal R$ with $|v|> (1-3\lambda)|r|$.

  The maximal arcs $p,e,f$ of $\Gamma$ subdivide $\beta$ as a
  concatenation
\[
\beta = p_1\dots p_s
\]
where $p_j$ are maximal arcs of $\Gamma$ (possibly traversed with the
opposite orientation) for $1<j<s$ and where $p_1,p_s$ are contained in
such maximal arcs.

There are two  cases to consider.

{\bf Case 1.} There is some $p_i$ such that $|p_i|\ge 6\lambda |r|$.

In this case $p_i$ is a subpath of $p^{\pm 1}$. After inverting
$\alpha$ and $w$ if needed, we may assume that in fact $p_i$ is a
subpath of $p$.  Note that the label on any part of $p_i$ which is
also read in $\beta$ in some different $p_j$ is a piece by definition.
Then the small cancellation condition $C'(\lambda)$
implies that there is a subsegment $q_i$ of $p$ such that $|q_i|\ge
3\lambda|r|$ and that $q_i$ does not overlap $p_j$ for $j\ne i$.

We then perform an $AO$-move by deleting the interior of the arc
$q_i$ and adding an arc labelled by the missing in $v$ part of
$r$, going from $o(\alpha)$ to $t(\alpha')$. This results in a
graph $\Gamma'$ with the smaller number of edges than in $\Gamma$.
By Corollary~\ref{tech}, the pair of elements of $G$ defined by
$\Gamma'$ is conjugate to the pair $(g_1,g_2)$. Note that
$\Gamma'$ is obtained from $\Gamma$ by removing a subsegment of
$p$ and then adding an arc connecting some vertex of one of the
two components of the remainder of $p$ to some vertex of the other
remaining component of $p$.

After removing the spikes ending in degree-one vertices if
necessary, we obtain another barbell-graph representing the
pair  $(g_1,g_2)$ but with a smaller number of edges than in
$\Gamma$ (where the label of the ``bar'' in this barbell need
not be reduced). This contradicts the minimal choice of $h$.

{\bf Case 2.} Suppose that $|p_i|< 6\lambda |r|$ for $1\le i\le s$.

Since $\alpha$ is cyclically reduced closed path that is not a single
loop-edge, it follows that for some $i$ we have $p_i=p^{\pm 1}$ and
hence $|p|<6\lambda |r|$.

Since $|v|> (1-4\lambda)|r|\ge |r|/2$, it follows that $v$ is
readable in a $u$-barbell graph with $|u|=|p|<6\lambda |r|$ and so
with $|u|/|v|< \frac{6\lambda}{2}\le \frac{1}{40}$.  This
contradicts the assumption that the tuple $(r_1,\dots, r_m)$
satisfies condition $Q_m(\lambda)$.

\end{proof}

  This main result now allows us to establish the desired rigidity
theorem.

\section{Rigidity of random quotients of the modular group}

\begin{conv}
For the remainder of this section, unless specified otherwise, let
$m\ge 1$ be an integer, let $0<\lambda \le \frac{1}{120}$ and let
\[
G=M/\langle\langle r_1,\dots, r_m \rangle\rangle, \tag{\dag}
\]
where the tuple $\tau=(r_1,\dots, r_m)$ of cyclically reduced
words satisfies the $Q_m(\lambda)$ condition.
\end{conv}

\begin{prop}\label{prop:basic} The following hold:
\begin{enumerate}
\item The group $G$ is one-ended and word-hyperbolic.
\item The cyclic subgroups $\langle a\rangle\le G$ or $\langle
  b\rangle\le G$ have orders $2$ and $3$, correspondingly, and these
  subgroups are malnormal in $G$.
\item Every nontrivial element of finite order in $G$ is conjugate to
  either $a$ or $b^{\pm 1}$.
\item Every finite subgroup of $G$ is conjugate to a subgroup of
  $\langle a\rangle$ or $\langle b\rangle$.
\item The center of $G$ is trivial.
\end{enumerate}
\end{prop}

\begin{proof}
  Statements (2),(3),(4) and (5) are straightforward applications of
  the small cancellation theory over free products, as explained in
  Section V.11 of Lyndon-Schupp ~\cite{LS}.  Also, Theorem~11.2 in Section~V.11 of
  \cite{LS} implies that presentation $(\dag)$ is a Dehn presentation
  and hence $G$ is word-hyperbolic.

  We already know from (2),(3),(4) that $G$ is not cyclic and hence
  has rank two. Suppose that $G$ is freely decomposable as $G=G_1\ast
  G_2$ where $G_i\ne 1$. Grushko's theorem then implies that each
  $G_i$ is 1-generated, that is, cyclic.  Since elements of finite
  order are always elliptic with respect to free product
  decompositions, it follows that $a$ is conjugate to an element of
  some factor $G_i$ and $b$ is conjugate to an element of some factor
  $G_j$. It follows that one of $G_1,G_2$ must be cyclic of order $2$
  and the other must be cyclic of order $3$. Hence
  $G=M/\langle\langle r_1,\dots,r_m  \rangle\rangle \simeq M$, which contradicts the fact that
  $M$ is Hopfian. Thus $G$ is freely indecomposable.

  Suppose now that $G$ admits a nontrivial splitting over a nontrivial
  finite group $H$. Since $G$ is generated by two elements of finite
  order, this splitting is not an HNN-extension. Thus $G=K\ast_H L$
  where $K,L\le G$ and $H\ne 1$, $H\ne K$, $H\ne L$.  Since $H\ne 1$ is finite,
  (4) implies that $H$ is conjugate to $\langle a\rangle$ or $\langle
  b\rangle$.  Therefore by (2) $H$ is malnormal in $G$. A theorem of
  Karrass and Solitar~\cite{KSo} then implies that $G$ cannot be
  generated by two elements, yielding a contradiction.

  Thus $G$ does not split nontrivially over a finite subgroup and
  hence, by Stallings' classic theorem~\cite{St}, $G$ is one-ended.
\end{proof}

Proposition~\ref{prop:main} and Proposition~\ref{prop:basic} immediately
imply:

\begin{thm}\label{thm:hom}
  Let $\psi: M\to G$ be a homomorphism. Then exactly one of the
  following mutually exclusive alternatives holds:
\begin{enumerate}
\item The map $\psi$ is injective but not onto.
\item The image of $\psi$ is a finite cyclic group of order $1,2$ or $3$.
\item The map $\psi$ is surjective and the pair $(\psi(a),\psi(b))$ is
  conjugate to $(a,b)$  or $(a,b^{-1})$ in $G$.
\end{enumerate}
\end{thm}

\begin{thm}\label{thm:endo} Let $\psi:G\to G$ be a homomorphism. Then the following hold:
\begin{enumerate}
\item If $\psi$ is onto then $\psi$ is injective. Hence $G$ is Hopfian.
\item If $\psi$ is injective then $\psi$ is onto. Hence $G$ is co-Hopfian.
\item If $\psi$ is an automorphism of $G$ then $\psi$ is inner. Hence
  $Out(G)=1$ and $G$ is a complete group.
\end{enumerate}
\end{thm}

\begin{proof}
Part (1)  follows directly from Theorem~\ref{thm:hom}. Indeed,
suppose $\psi: G\to G$ is an onto endomorphism. Then by
Theorem~\ref{thm:hom} we know that $(\psi(a),\psi(b))$ is
  conjugate to $(a,b)$  or $(a,b^{-1})$ in $G$. Suppose the latter
  holds. Then $\eta(r_1)=_G 1$. Hence by the small cancellation
  assumption on $\tau$ it follows that $\eta(r_1)$ contains a
  subword that is more than a half of a cyclic permutation of some
  $r_i^{\pm 1}$. This contradicts condition (5) in Definition~\ref{defn:Q} of
  $Q_m(\lambda)$. Thus $(\psi(a),\psi(b))$ is
  conjugate to $(a,b)$ in $G$. Therefore $\psi$ is an inner
  automorphism of $G$ and in particular, $\psi$ is injective.

For part (2), suppose that $\psi:G\to G$ is an injective
endomorphism. We need to show that $\psi$ is onto.

Theorem~\ref{thm:hom} implies that either the image of $\psi$ is
isomorphic to $M$ or, after a post-composition with an inner
automorphism we have $(\psi(a),\psi(b))=(a,b^{\pm 1})$. The former
is impossible since $G$ is one-ended while $M$ is not. In the
latter case the image of $\psi$ is generated by a tuple conjugate
to $(a,b^{\pm 1})$ and hence $\psi$ is onto, as required.

For part (3), let $\psi: G\to G$ be an automorphism of $G$. Since
$G$ is one-ended, Proposition~\ref{prop:main} implies that
$(\psi(a), \psi(b))$ is conjugate in $G$ to either $(a,b)$ or
$(a,b^{\-1})$. In the former case $\psi$ is inner, as required. In
the latter case we obtain a contradiction, exactly as in the proof
of part (1).
\end{proof}

  The following result is essentially due to Greendlinger~\cite{Green} who
  proved it in the context of small cancellation quotients of a free
  group.  We present an argument for completeness.

\begin{prop} [Greendlinger's Theorem] \label{prop:green}
  Let $R_1,R_2$ be finite nonempty symmetrized sets of cyclically
  reduced words in $M$ such that each $R_i$ satisfies the $C'(1/8)$
  small cancellation condition and such that 
$\langle\langle R_1\rangle\rangle=\langle\langle R_2 \rangle\rangle$ 
 in   $M$.   Then $R_1=R_2$.
\end{prop}

\begin{proof}

  Suppose that the result fails and that $R_1\ne R_2$. Let $r$ be the
  shortest element from the symmetric difference of $R_1$ and $R_2$.
  Without loss of generality we may assume that $r\in R_1-R_2$.

  Since $r\in \langle \langle R_2 \rangle\rangle$, the normal closure of $R_2$.
  there exists a reduced van Kampen diagram $\Delta$
  over $R_2$ with $r$ being the label of the boundary cycle. If $\Delta$
  contains a single region then $r\in R_2$, contrary to our
  assumptions.  Thus $\Delta$ contains at least two regions. Since $R_2$
  satisfies the $C'(1/8)$-small cancellation condition, the perimeter
  of $\Delta$ is longer than the perimeter of every region in $\Delta$. Thus all
  regions of $\Delta$ have boundaries labelled by elements of $R_2$ that
  are shorter than $r$. The minimal choice of $r$ implies that these
  elements of $R_2$ also belong to $R_1$.

  Thus $\Delta$ is a reduced diagram over $R_1$ with the boundary cycle $r$
  and with at least two regions. Again since $R_1$ is $C'(1/8)$, it
  follows that there is a region $D_0$ of $\Delta$ with boundary cycle
  labelled by $r'\in R_1$ such that $|r'|<|\partial D|=|r|$ and such
  that there is an arc in the boundary cycle of $D_0$ that is
  contained in the boundary cycle of $\Delta$ and such that this arc has
  length at least $|r'|/2$. Since $r,r'\in R_1$ and $|r'|<|r|$, this
  contradicts the $C'(1/8)$-small cancellation condition for $R_1$.
\end{proof}

\begin{thm}\label{thm:rigid}[Isomorphism Rigidity for random quotients]
  Let $m\ge 1$.  Then for any $\tau=(r_1,\dots, r_m)$ satisfying
  condition $Q_m(\frac{1}{120})$ and any finite symmetrized
  $S\subseteq \mathcal C$ satisfying the standard small cancellation
  condition $C'(1/8)$ we have $G_\tau\simeq   M/\langle\langle S \rangle\rangle$
  if and only if $R(\tau)=S$ or $R(\tau)=\eta(S)$.
\end{thm}
\begin{proof}
  The ``if'' direction is obvious.  Suppose now that $G_\tau\simeq
  M/\langle\langle S \rangle\rangle$ and let $\psi:M/\langle\langle S \rangle\rangle\to G_\tau$ be an isomorphism.
  Since $\psi$ is onto, Theorem~\ref{thm:hom} implies that the pair
  $(\psi(a),\psi(b))$ is conjugate to $(a,b)$ or $(a,b^{-1})$ in
  $G_\tau$. After composing $\psi$ with an inner automorphism of
  $G_\tau$, we may assume that $\psi(a)=a$ and $\psi(b)=b^{\pm 1}$.
  Then in $M$ we have
\[
\langle\langle r_1,\dots, r_m \rangle\rangle=\langle\langle
\eta^{\epsilon}(S) \rangle\rangle
\]
for some $\epsilon\in\{0,1\}$. Since both $R(\tau)$ and $S$ satisfy
the $C'(1/8)$ small cancellation condition over $M$,
Proposition~\ref{prop:green} implies that
$R(\tau)=\eta^{\epsilon}(S)$, as required.
\end{proof}

Theorem~\ref{thm:rigid} and the definition of $Q_m(\lambda)$
immediately imply:

\begin{cor}\label{thm:rigid1}
  Let $m\ge 1$ and let   $\sigma=(r_1,\dots, r_m), \tau =(s_1,\dots, s_m)
   \in Q_m ( \frac{1}{120})$.
   Then $G_\sigma\simeq G_\tau $ if and only if there exists a
  reordering $\tau'=(s_1',\dots, s_m')$ of $\tau$ and $\delta \in\{0,1\}$ such that each
  $s_i'$ is a cyclic permutation of $\eta^\delta(r_i)$ or $\eta^\delta(r_i^{-1})$ for
  $i=1,\dots, m$.
\end{cor}

\begin{cor}\label{cor:hom}
Let $m\ge 1$ be an integer and let $\sigma, \tau \in
U_m(\frac{1}{120})$ be such that $|\sigma|\le |\tau|$. Suppose
that $\psi:G_\sigma \to G_\tau$ is a homomorphism. Then the
following hold:

\begin{enumerate}
\item If  $|\sigma|<|\tau$ then $\psi(G_\sigma)$ is a finite cyclic
group of order at most $3$.

\item If $|\sigma|=|\tau|$ then either $\psi(G_\sigma)$ is a
finite cyclic group of order at most $3$ or $\psi$ is an
isomorphism. In the latter case, after a possible post-composition
of $\psi$ with an inner automorphism of $G_\tau$, we have
$\psi(a)=\eta^\delta(a)=a$ and $\psi(b)=\eta^\delta(b)=b^{\pm 1}$
for some $\delta\in \{0,1\}$ and, moreover,
$\eta^\delta(R(\sigma))=R(\tau)$.
\end{enumerate}

\end{cor}

\begin{proof}

(1) Suppose that $|\tau|<|\sigma|$. Since $M$ is Hopfian,
Theorem~\ref{thm:hom} implies that either $\psi(G_\tau)$ is a
finite cyclic group of order at most $3$ or, after composing
$\psi$ with an inner automorphism we have $\psi(a)=a$ and
$\psi(b)=b^{\pm 1}$.  Hence for some $\delta\in \{0,1\}$ we have
$\eta^\delta(\tau)\subseteq \langle\langle \sigma \rangle\rangle$.
However, $R(\sigma)$ satisfies the $C'(1/8)$ small cancellation
condition, and therefore every nontrivial element from
$\langle\langle \sigma \rangle\rangle$ in $M$ has length $\ge
|\sigma|$. On the other hand, $\eta^\delta(\tau)\subseteq
\langle\langle \sigma \rangle\rangle$ and $|\tau|<|\sigma|$,
yielding a contradiction. This proves part (1) of
Corollary~\ref{cor:hom}.

(2) Suppose now that $|\tau|=|\sigma|$. Again, since $M$ is
Hopfian, Theorem~\ref{thm:hom} implies that either $\psi(G_\tau)$
is a finite cyclic group of order at most $3$ or, after composing
$\psi$ with an inner automorphism we have $\psi(a)=a$ and
$\psi(b)=b^{\pm 1}$. Thus, as before for some $\delta\in \{0,1\}$
we have $\eta^\delta(\tau)\subseteq \langle\langle \sigma
\rangle\rangle$. Since $R(\sigma)$ satisfies the $C'(1/8)$ small
cancellation condition, if a nontrivial cyclically reduced word
$r$ belongs to $\langle\langle \sigma \rangle\rangle$ then either
$|r|>|\sigma|$ or $r\in R(\sigma)$. Since $|\tau|=|\sigma|$ and
$\eta^\delta(\tau)\subseteq \langle\langle \sigma \rangle\rangle$,
it follows that $\eta^\delta(\tau)\subseteq R(\sigma)$ and hence
$\eta^\delta(R(\tau))\subseteq R(\sigma)$. By definition of
$U(\lambda)$ it follows that
\[
\#R(\sigma)=\#R(\tau)=\#\eta^\delta(R(\tau))= 2mn
\]
where $n=|\tau|=|\sigma|$. Hence $\eta^\delta(R(\tau))= R(\sigma)$
and therefore $\psi$ is an isomorphism, as claimed.
\end{proof}

\begin{cor}\label{thm:not}
  Let $p>m\ge 1$ be integers. Let $\sigma=(r_1,\dots, r_m)\in
  U_m(\frac{1}{120})$ and $\tau=(s_1,\dots, s_p)\in
  U_p(\frac{1}{120})$ be such that $|\tau|=|\sigma|$. Then
  $G_\tau\not\simeq G_\sigma$.
\end{cor}
\begin{proof}
  Let $n:=|r_i|=|s_j|$.

  Suppose that $G_\tau\simeq G_\sigma$ and let $\psi:G_\sigma\to
  G_\tau$ be an isomorphism. Then by Theorem~\ref{thm:rigid} either
  $R(\tau)=R(\sigma)$ or $R(\tau)=\eta(R(\sigma))$. The definition of
  $U_m(\frac{1}{120}$ and $U_p(\frac{1}{120}$ implies that
\[
\#R(\tau)= 2mn, \qquad \#R(\sigma)=\#\eta(R(\sigma))= 2pn
\]
This is a contradiction since $m\ne p$ and thus $2pn \ne 2mn$.
\end{proof}

  The following sections of the paper can each be read in any order.
Both essential incompressibility and counting isomorphism types
 follow directly from rigidity  and are independent of each other.

\section{Counting the isomorphism types of random quotients}

\begin{defn}
  For an even integer $n=2k>0$ let $I_m(n)$ denote the number of
  isomorphism types of groups given by presentations
\[
M/\langle \langle r_1,\dots,r_m \rangle\rangle,
\]
where $r_1,\dots, r_m$ are cyclically reduced words of length $n$ in
$M$
\end{defn}

\begin{thm}~\label{thm:types}
  We have for even $n\to\infty$
\[
I_m(n)\sim \frac{(2^{n/2+1})^{m}}{2\ m!(2n)^m},
\]
that is
\[
\lim_{k\to\infty}\frac{2 I_{m}(2k)\ m!(4k)^m}{(2^{k+1})^{m}}=1.
\]
\end{thm}
\begin{proof}
  Let $V_m:=U_m(\frac{1}{120}) \cap \mathcal T_m$. Then $V_m$ is
  exponentially generic in $\mathcal T_m$.  Let $X_m$ be
  the complement of $V_m$ in $\mathcal T_m$.

  Recall that  the number of cyclically reduced words of even length
  $n = 2k$ is $2~2^k = 2^{k+1}$.
  For every $\tau=(r_1,\dots, r_m)\in V_m$ with $|r_i|=n$
  there are precisely $2\ m!(2n)^m$ tuples $\sigma\in Q_m'$ such that
  $G_\tau\simeq G_\sigma$. Let $\gamma(n,V_m)$ be the number of
  tuples from $V_m$ with entries of length $n$. Then the tuples
  $\tau$ from $V_m$ with entries of length $n$ represent exactly
  $I_m'(n):=\frac{\gamma(n,V_m)}{2\ m!(2n)^m}$ distinct isomorphism
  types of groups $G_\tau$. Let $I_m''(n)$ be the number of
  isomorphism types of groups $G_\tau$ corresponding to $\tau\in X_m$
  with entries of length $n$. Thus $I_m(n)=I_m'(n)+I_m''(n)$.

Since $V_m$ is exponentially generic in $\mathcal T_m$,
we have
\[
\lim_{k\to\infty}\frac{\gamma(2k,V_m)}{(2^{k+1})^{m}}=1
\]
and
\[
\lim_{k\to\infty}\frac{\gamma(2k,X_m)}{(2^{k+1})^{m}}=0\tag{*}
\]
with exponentially fast convergence.
Therefore
\[
\lim_{k\to\infty}\frac{2 I_{m}'(2k)\  m!(4k)^m}{(2^{k+1})^{m}}=1.
\]
Moreover $I_m''(2k)\le \gamma(2k,X_m)$ and, since the convergence in
$(*)$ is exponentially fast, we have
\[
\lim_{k\to\infty}\frac{2 I_{m}''(2k)\  m!(4k)^m}{(2^{k+1})^{m}}=0.
\]
Since $I_m(n)=I_m'(n)+I_m''(n)$, we conclude that
\[
\lim_{k\to\infty}\frac{2 I_{m}(2k)\  m!(4k)^m}{(2^{k+1})^{m}}=1,
\]
as required.
\end{proof}

\section{Kolmogorov Complexity and the $T$-invariant}

Intuitively, the Kolmogorov complexity $K(x)$ of a finite binary
string $x$ is the size of the smallest computer program $M$ that
can compute $x$. In order to make this notion precise  one needs
to first fix a ``programming language'' but it turns out that all
reasonable choices yield measures which are equivalent up to an
additive constant. We refer the reader to the book of Li and
Vitanyi~\cite{LV} for a detailed treatment of the subject and we
recall  only a few relevant facts and definitions here.

\begin{defn}\label{defn:kolm}
  Fix a universal Turing machine $U$ with the alphabet \newline
  $\Sigma:=\{0,1\}$. Then $U$ computes a universal partial recursive
  function $\phi$ from $\Sigma^{\ast}$ to $\Sigma^{\ast}$. That is,
  for any partial recursive function $\psi$ there is a string $z\in
  \Sigma^{\ast}$ such that for all $x\in \Sigma^{\ast}$, $\phi (zx) =
  \psi(x)$.

  For a finite binary string $x\in \Sigma^{\ast}$ we define the
  \emph{Kolmogorov complexity} $K(x)$ as

\[
K(x):=\min \{ |p| : p\in \Sigma^{\ast}, \phi(p)=x\}.
\]

Let $B$ be another finite nonempty alphabet.  We choose a recursive
bijection $h: B^{\ast}\to \{0,1\}^{\ast}$.

For any string $x\in B^*$ define its \emph{Kolmogorov complexity}
$K_B(x)$ as

\[
K_B(x):=K(h(x)).
\]

Let $\widehat A$ denote the alphabet consisting of $A$ together with the
extra symbols ``('', ``)'' and ``,''. We regard
tuples $\tau=(r_1,\dots, r_m)\in \mathcal C^m$ as words in the
alphabet $\widehat A$. For an $m$-tuple $\tau=(r_1,\dots, r_m)\in \mathcal
C^m$ define the \emph{Kolmogorov complexity $K(\tau)$} as
\[
K(\tau):=K_{\widehat A}(\tau).
\]

\end{defn}

\begin{rem}
The notations in the present paper differ slightly from those use
in \cite{KSn}. In \cite{KSn} $K(x)$ stood for the prefix complexity of
$x$ while the Kolmogorov complexity of $x$ was denoted by $C(x)$.
\end{rem}

\begin{rem}[The General Enumeration Argument]

  Let $\mathcal P$ be a recursively enumerable class of finite group
  presentations. Then there is an algorithm which, when given a finite
  presentation of a group $G$ known to be isomorphic to a
  group defined by some  presentation from $\mathcal P$, actually finds
  a presentation $\Pi$ of $G$ which is in $\mathcal P$.

  Say $G$ is given by a presentation $G=\langle X|R\rangle$. We start
  enumerating all presentations $\Pi=\langle Y| S\rangle\in \mathcal
  P$. For each such presentation we start enumerating pairs of maps
  $(\alpha,\beta)$ where $\alpha:X\to F(Y)$ and $\beta:Y\to F(X)$. For
  each such pair start enumerating $\langle\langle S \rangle\rangle\le F(Y)$ and $\langle\langle R \rangle\rangle\le
  F(X)$ and start checking if $\alpha, \beta$ extend to homomorphisms
  $\alpha: G\to G(\Pi)$ and $\beta: G(\Pi)\to G$, where $G(\Pi)$ is the group defined by $\Pi$.
  If we find a pair
  $(\alpha,\beta)$ where both $\alpha,\beta$ extend to such
  homomorphisms, we start checking if $x^{-1}\beta\alpha(x)\in \langle\langle R \rangle\rangle$
  and $y^{-1}\alpha\beta(y)\in \langle\langle S \rangle\rangle$ for all $x\in X,y\in Y$. If
  yes, then $\alpha: G\to G(\Pi)$ and $\beta: G(\Pi)\to G$ are
  mutually inverse isomorphisms. Since $G$ is known to be isomorphic
  to some group defined by a presentation from $\mathcal P$, this
  process is guaranteed to terminate.

  We  refer to the above algorithm constructed above as the \emph{general
    enumeration algorithm} for $\mathcal P$.
\end{rem}

We will need the following statement that follows directly from
Proposition~2.5 and Lemma~2.7 in \cite{KSn}.

\begin{prop}\label{prop:old}
  Let $m\ge 1$ be a fixed integer. Let $\Omega\subseteq \mathcal
  T_m$ be a nonempty subset equipped with a discrete non-vanishing
  probability measure $P$, so that $\sum_{\tau\in \Omega}
  P(\{\tau\})=1$.  Denote $\mu(\tau):=P(\{\tau\})$ for any $\tau\in
  \Omega$.

  Then for any $\delta>0$ we have

\[
P\big(2K(\tau)\ge -\log_2 \mu(\tau) -\log_2 \delta -c \big)\ge
1-\frac{1}{\delta}.
\]
where $c=c(m)>0$ is a constant independent of $\Omega, P$.
\end{prop}

\begin{cor}\label{cor:main}
  Let $\delta>0$ and let $Z_\delta\subseteq \mathcal T_m$ be the set
  of all tuples $\tau \in \mathcal T_m$ such that
\[
2K(\tau)\ge m(n/2+1)-\log_2 \delta -c,
\]
where $n$ is the length of each entry of $\tau$.  Then for all even
$n>0$ we have
\[
\frac{\gamma(n,Z_\delta)}{\gamma(n,\mathcal T_m)}\ge
1-\frac{1}{\delta}.
\]
\end{cor}
\begin{proof}
  Let $n>0$ be an even integer and let $\Omega$ be the set of all
  tuples in $T_m$ with entries of length $n$. Let $P$ be the
  uniform probability measure on $\Omega$. Then for every $\tau\in
  \Omega$ we have
\[
\mu(\tau)=\frac{1}{2^{m(n/2+1)}}.
\]
Hence
\[
-\log_2 \mu(\tau)=m(n/2+1).
\]
Applying Proposition~\ref{prop:old} we get
\[
P\left( 2K(\tau)\ge m(n/2+1)-\log_2 \delta -c\right)\ge
1-\frac{1}{\delta},
\]
as required.
\end{proof}

\begin{defn}
  Let $0<\lambda\le \frac{1}{120}$ and let $m\ge 1$ be an integer.
  Let $U_m'(\lambda)$ denote the set of all $\tau=(r_1,\dots, r_m)\in
  U_m(\lambda)$ with the following property.

  Whenever $\tau'=(r_1',\dots, r_m')$ is obtained from $\tau$ by a
  formally nontrivial combination of reordering the elements of
  $\tau$, taking cyclic permutations and possible inverses of its
  entries and possibly applying $\eta$ to the tuple, then
\[
(v_1,\dots, v_m)\ne (v_1',\dots, v_m')
\]
where $v_i$ is the initial segment of $r_i$ of length $\lfloor \lambda
n\rfloor$, where $v_i'$ is the initial segment of $r_i'$ of length
$\lfloor \lambda n\rfloor$ and where $n=|r_i|=|r_j'|$.
\end{defn}

\begin{prop}\label{prop:unique}
  Let $0<\lambda\le \frac{1}{120}$ and let $m\ge 1$ be an integer.
  Then $U_m'(\lambda)$ is exponentially generic in
  $\mathcal T_m$.
\end{prop}
\begin{proof}
  The proof is a straightforward generalization of the proofs of
  Lemmas 4.5, 4.6 and 4.8 in \cite{KSn}, where a similar statement was
  considered for the case of a free group. We leave the details to the
  reader.
\end{proof}

\begin{lem}\label{lem:N}
  There exists a constant $N=N(m)>0$ with the following property.  Let
  $0<\lambda \le 1/24$ be a rational number and let $\tau\in
  U_m'(\lambda)$ be a tuple consisting of words of length $n\ge 2$.

  Suppose $G_\tau$ can be presented by a finite presentation
\[
\Pi=\langle b_1, \dots b_s | w_1, \dots, w_t\rangle \tag{\ddag}
\]
where $t\ge 1$.

Then $K(\tau)\le N\ell_1(\Pi)\log_2 \ell_1(\Pi)+nN\lambda+N$.
\end{lem}

\begin{proof}
  We describe an algorithm $\mathcal A$, which, given a presentation
  $(\ddag)$ for $G_\tau$ and a tuple $(u_1,\dots, u_m)$ of initial
  segment $u_i$ of $r_i$ of length $\lfloor\lambda n\rfloor$, will
  recover the tuple $\tau$.

  First, note that we are assuming that $(\ddag)$ defines a group
  isomorphic to a group $G_\sigma$ for some $\sigma\in U_m'(\lambda)$.
  Since the set $U_m'(\lambda)$ is recursive, by a general enumeration
  algorithm we can algorithmically find some $\sigma\in U_m'(\lambda)$
  such that $(\ddag)$ defines a group isomorphic to $G_\sigma$ and
  hence to $G_\tau$.

  We then perform all possible ways of applying to $\sigma$ a
  combination of reordering of the tuple entries, cyclic permutations
  and possible inversions of its entries and possibly applying $\eta$
  to the tuple. Corollary~\ref{thm:rigid1} implies that this
  collection of tuples will contain $\tau$.

  For each of the resulting tuples we record the sequence of initial
  segments of its entries of length $\lfloor\lambda n\rfloor$ and
  compare it with $(u_1,\dots, u_m)$. The definition of
  $U_m'(\lambda)$ that there will be exactly one tuple for which this
  sequence of initial entries coincides with $(u_1,\dots, u_m)$ and
  this tuple is $\tau$.

  The general enumeration algorithm is fixed.  The further input of
  $\mathcal A$, required to compute $\tau$, consists of the
  presentation $(\ddag)$ and tuple of the initial segments
  $(u_1,\dots, u_m)$ of $r_1,\dots, r_m$ with $|u_i|=\lfloor \lambda n
  \rfloor$.

  We want to estimate the length of this input when expressed as a
  binary sequence. Put $T=\ell_1(\Pi)$. First note that in $(\ddag)$
  every $b_i$ must occur in some $w_j^{\pm 1}$ since $G_\tau$ is a
  one-ended group by Theorem~\ref{thm:rigid} and therefore $s\le T$.

  We can now encode the presentation $(\ddag)$ by writing each
  subscript $i=1,\dots, s$ for each occurrence of $b_i$ in $(\ddag)$
  as a binary integer. Using $\overline i$ to denote the binary
  expression for $i$, we replace each occurrence of $b_i$ in $(\ddag)$
  by $b \overline i$ and each occurrence of $b_i^{-1}$ by $-b\overline
  i$.  Note that the bit-length of the binary expression $\overline i$
  of $i$ is at most $\log_2 i$.  This produces an unambiguous encoding
  of $(\ddag)$ as a string $W$ of length at most $O(T\log_2 T)$ over
  the six letter alphabet
  \[
  b\quad 0\quad 1\quad -\quad , \quad |
  \] and this
  alphabet can then be block-coded into binary in the standard way.

  Since the number $m$ and the alphabet $A=\{a,b,b^{-1}\}$ are fixed,
  describing $(u_1,\dots,u_m)$ requires at most $O(\lambda n)$ number
  of bits.

  Hence there exist a constant $N=N(m)>0$ independent of $\tau$ such
  that
\[
K(\tau)\le NT\log_2 T+Nn\lambda+N.
\]

\end{proof}

\begin{thm}\label{thm:kolm}
  Let $m\ge 1$ be a fixed integer.  For any $0<\epsilon<1$ there is an
  integer $n_0>0$ and a constant $L=L(m,\epsilon)>0$ with the
  following property.

  Let $J$ be the set of all tuples $\tau\in \mathcal T_m$ such that
\[
T_1(G_\tau)\log_2 T_1(G_\tau) \ge L |\tau|.
\]
Then for any $n\ge n_0$

\[
\frac{\gamma(n,J)}{\gamma(n,\mathcal T_m)}\ge 1-\epsilon.
\]
\end{thm}

\begin{proof}

  Let $N>0$ be the constant provided by Lemma~\ref{lem:N}. Choose a
  rational number $\lambda$, $0<\lambda<2/135$ so that $\frac{m}{4}-
  N\lambda>0$.

  Let $0<\epsilon<1$ be arbitrary and let $\delta>0$ be such that
  $\frac{2}{\delta}<\epsilon$.

  As in Corollary~\ref{cor:main} let $Z_\delta$ be let
  $Z_\delta\subseteq \mathcal T_m$ be the set of all tuples $\tau
  \in \mathcal T_m$ such that
\[
K(\tau)\ge \frac{1}{2}[m(n/2+1)-\log_2 \delta -c],
\]
where $n$ is the length of each entry of $\tau$.  Then for all even
$n>0$ we have
\[
\frac{\gamma(n,Z_\delta)}{\gamma(n,\mathcal T_m)}\ge
1-\frac{1}{\delta}.
\]

Since by Proposition~\ref{prop:unique} $U_m'(\lambda)$ is
exponentially generic in $\mathcal T_m$, there is $n_1>0$ such that
for any $n\ge n_1$

\[
\frac{\gamma(n,Z_\delta\cap U_m'(\lambda))}{\gamma(n,\mathcal
  T_m)}\ge 1-\frac{2}{\delta}\ge 1-\epsilon.
\]

Now suppose $\tau\in Z_\delta\cap U_m'(\lambda)$ and
$n:=|\tau|\ge n_1$.

Then by Lemma~\ref{lem:N}
\[
\frac{1}{2}[m(n/2+1)-\log_2 \delta -c]\le K(\tau)\le
NT_1(G_\tau)\log_2 T_1(G_\tau)+Nn\lambda+N.
\]

Hence for $n\ge n_1$
\[
(\frac{m}{4}-N\lambda)n +\frac{m}{2}-\log_2 \delta -c-N\le
NT_1(G_\tau)\log_2 T_1(G_\tau).
\]

This implies the statement of the theorem.

\end{proof}

\end{document}